\documentclass[titlepage,twoside,12pt]{article}
\usepackage{amssymb}
\usepackage{amsfonts}
\usepackage{amsmath}
\begin{document}

{\bf \Large BRIEF NOTES ON SHEAVES THEORY} \\ \\

{\bf Elemer E Rosinger} \\
Department of Mathematics \\
and Applied Mathematics \\
University of Pretoria \\
Pretoria \\
0002 South Africa \\
eerosinger@hotmail.com \\ \\

{\bf Abstract} \\

Essentials of sheaves are briefly presented, followed by related comments on presheaves, bundles, manifolds and
singularities, aiming to point to their differences not only in their different formal mathematical structures, but also
in the very purposes for which they were introduced in the first place. \\ \\ \\

{\bf 1. DEFINITIONS} \\

{\bf Definition 1.1. SHEAF} \\

${\cal S} = ( {\cal S}, \pi, X )$ is a {\it sheaf over} $X$, where $X$ and ${\cal S}$ are
topological spaces, if and only if \\

(1.1)~~~ $ \pi : {\cal S} \longrightarrow X $ \\

is a surjective local homeomorphism. \\

{\bf Lemma 1.1.} \\

$ \left \{~ B \subseteq {\cal S} ~~~~ \begin{array}{|l}
                                   ~~1)~~ B ~~\mbox{open} \\
                                   ~~2)~~ \pi|_B : B \longrightarrow \pi ( B ) \subseteq X
                                            ~~\mbox{is a local homeomorphism}
                                 \end{array} ~ \right \} $ \\

is a basis of the topology on ${\cal S}$. \\

{\bf Definition 1.1a. FIBERS} \\

$ {\cal S}_x = \pi^{-1} ( \{ x \} ) $ \\

is the {\it fiber} over $x \in X$.

\hfill $\Box$ \\

We have \\

$ {\cal S} = \sum_{x \in X} {\cal S}_x $ ~ a union of pair-wise disjoint sets \\

$ {\cal S}_x ~~\mbox{discrete subset of}~ {\cal S},~~~ x \in X $ \\

{\bf Example 1. CONSTANT SHEAF} \\

Let $M$ be any set considered with the discrete topology. Let ${\cal S} = X \times M$ and
$\pi = pr_X : {\cal S} \longrightarrow X$. Then $( {\cal S}, \pi, X )$ is a sheaf on $X$ and
${\cal S}_x = \{ x \} \times M$, for $x \in X$. \\

{\bf Remark 1.1. Difference between sheaves and fiber bundles} \\

A {\it bundle} is given by $( E, \pi, B )$, where $E$ and $B$ are topological spaces and $\pi
: E \longrightarrow B$ is a surjective continuous mapping. The case of interest is when there
exists a topological space $F$, called {\it fiber}, such that \\

$ \pi^{-1} ( x ) $ ~is homeomorphic to~ $ F $, ~for~ $ x \in B $ \\

The main property of a bundle is its {\it local homeomorphism} with a {\it Cartesian
product}, namely \\

$ \begin{array}{l}
          \forall~~ z \in E ~:~ \\ \\
          \exists~~ x = \pi ( z ) \in U \subseteq B,~ U ~~\mbox{open} ~:~ \\ \\
          ~~~~ \pi^{-1} ( U ) \subseteq E ~~\mbox{is homeomorphic with}~ U \times F
   \end{array} $ \\ \\

Thus unlike in a sheaf $( {\cal S}, \pi, X )$, where ${\cal S}$ is locally homeomorphic with
$X$, in a bundle $( E, \pi, B )$, the space $E$ is {\it not} locally homeomorphic with $B$,
but with $U \times F$, for suitable open subsets $U \subseteq B$. \\

{\bf Definition 1.2. SHEAF MORPHISM} \\

Let $( {\cal S}, \pi, X ),~ ( {\cal S\,'}, \pi\,', X )$ be two sheaves on $X$. A mapping
$\phi : {\cal S} \longrightarrow {\cal S\,'}$ is a {\it sheaf morphism}, if and only if the
diagram commutes \\

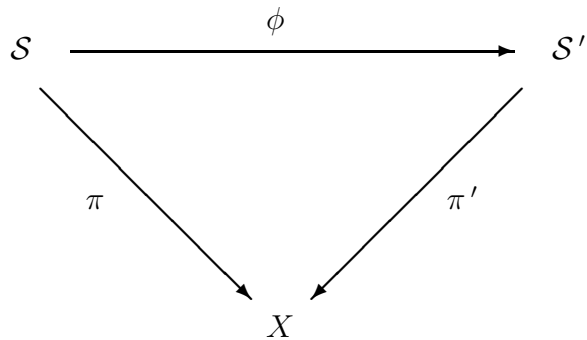
\begin{math}
\setlength{\unitlength}{0.2cm}
\thicklines
\begin{picture}(60,20)

\put(11,16){${\cal S}$}
\put(28,18){$\phi$}
\put(15,16.5){\vector(1,0){29.5}}
\put(47,16){${\cal S\,'}$}
\put(16,6){$\pi$}
\put(13,14){\vector(1,-1){14}}
\put(28,-2.5){$X$}
\put(40,6){$\pi\,'$}
\put(45,14){\vector(-1,-1){14}}

\end{picture}
\end{math} \\ \\

We denote by \\

$ {\cal S}h_X $ \\

the category of sheaves over $X$.

\hfill $\Box$ \\

We have \\

$ \phi ~~\mbox{sheaf morphism} ~~\Longleftrightarrow~~ (~ \phi ( {\cal S}_x ) \subseteq
                         {\cal S}\,'_x,~~~ x \in X ~) $ \\

{\bf Lemma 1.2.} \\

If $\phi$ is sheaf morphism, then it is a local homeomorphism from ${\cal S}$ to
${\cal S\,'}$. \\ \\

{\bf 2. SHEAF SECTIONS} \\

{\bf Definition 2.1. SECTIONS} \\

A {\it section} in the sheaf $( {\cal S}, \pi, X )$ is any continuous mapping $s : U
\longrightarrow {\cal S}$, where $U \subseteq X$ is open, and the diagram commutes

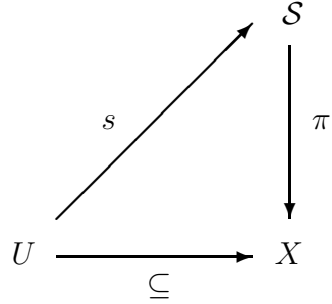
\begin{math}
\setlength{\unitlength}{0.2cm}
\thicklines
\begin{picture}(60,20)

\put(27,16){${\cal S}$}
\put(29,9){$\pi$}
\put(27.5,14.5){\vector(0,-1){11.5}}
\put(26.5,0){$X$}
\put(15,9){$s$}
\put(12,3){\vector(1,1){13}}
\put(9,0){$U$}
\put(12,0.5){\vector(1,0){13}}
\put(18,-2){$\subseteq$}

\end{picture}
\end{math} \\

thus we have \\

$ \pi \circ s = id_U $ \\

We denote \\

$ {\cal S} ( U ) = \Gamma ( U, {\cal S} ) $ \\

the set of all such sections. \\

{\bf Lemma 2.1.} \\

If $s \in \Gamma ( U, {\cal S} )$, then \\

1)~~ $ s ( U ) $ ~is open in~ $X $ \\

2)~~ $ s : U \longrightarrow s ( U ) \subseteq {\cal S} $ ~is a homeomorphism \\

3)~~ if~ $ V \subseteq {\cal S} $ ~is open and~ $ \pi|_V : V \longrightarrow \pi ( V )
\subseteq X $ ~is a local

~~~~~~homomorphism, then \\

$~~~~~~ s = ( \pi|_V )^{-1} : \pi ( V ) \longrightarrow V $ ~is a section \\

4)~~ $ \forall~~ z \in {\cal S} ~:~ \exists~~ x = \pi ( z ) \in U \subseteq X $ ~open,~
$ s \in \Gamma ( U, {\cal S} ) ~:~ z = s ( x ) $ \\

5)~~ if~ $ s \in \Gamma ( U, {\cal S} ),~ t \in \Gamma ( V, {\cal S} )$ ~and for some~
$x \in U \cap V$ ~we

~~~~~~have~ $ s ( x ) = t ( x ) $, ~then there exists~ $ x \in W \subseteq U \cap V $ ~open,
such

~~~~~~that~ $ s_W = t_W $. \\

6)~~ if~ $ A \subseteq {\cal S}$ ~is open, such that~ $ \pi|_A $ ~is a homeomorphism, then

~~~~~~~there exists~ $ U \subseteq X $ ~open, and~ $ s \in \Gamma ( U, {\cal S} ) $, ~such that~
$ A = s ( U ) $. \\

7)~~ $ \{~ s ( U ) ~|~ U \subseteq X ~~\mbox{open}, ~ s \in \Gamma ( U, {\cal S} ) ~\} $ ~is a
                   basis for the topology

~~~~~~of ${\cal S}$. \\

{\bf Remark 2.1.} \\

1) From 2) above follows that sheaves have {\it plenty} of {\it sections}. \\

2) As seen in section 3, a sheaf can in fact be {\it determined} by its {\it sections}. \\

3) From 4) and 5) above follows that every $ z \in {\cal S}$ is the {\it germ} at $x =
\pi ( z ) \in X$ of a suitable section $s$. \\

4) Thus a sheaf ${\cal S}$ is {\it determined} by the {\it germs} of its sections. \\

5) In fact, as seen in 7) above, the {\it topology} of ${\cal S}$ is also {\it determined} by
the {\it sections}. \\

{\bf Proposition 2.1.} \\

Let $( {\cal S}, \pi, X ),~ ( {\cal S\,'}, \pi\,', X )$ be two sheaves on $X$, and a
continuous mapping $\phi : {\cal S} \longrightarrow {\cal S\,'}$. Then the following are
equivalent \\

1)~~ $ \phi $ ~is a sheaf morphism \\

2)~~ for every open~ $ U \subseteq X $ ~we have~ $ \Gamma ( U, {\cal S} ) \ni s \longmapsto
           \phi \circ s \in \Gamma ( U, {\cal S}\,' ) $ \\

3)~~ for every~ $ z \in {\cal S} $ ~there exists an open~ $ U \subseteq X $ ~and~ $ s \in
           \Gamma ( U, {\cal S} ) $,

~~~~~~such that~ $ z \in s ( U ),~ \phi \circ s \in \Gamma ( U, {\cal S}\,' ) $ \\ \\

{\bf 3. A SHEAF IS ITS SECTIONS} \\ \\

$\begin{array}{||l}
          ~~~"FIBERS" + "SECTIONS" ~~------>~~ SHEAF
  \end{array} $ \\ \\

Let $X$ be a topological space, and \\

(3.1)~~~ $ S_x,~~~ x \in X $ \\

a pair-wise disjoint family of sets, while \\

(3.2)~~~ $ \sigma =  \{~ s = ( U, s, {\cal S} ) ~|~ U \subseteq X,~ U ~~\mbox{open},~ s : U
                               \longrightarrow {\cal S} ~\} $ \\

where \\

(3.3)~~~ $ {\cal S} = \sum_{x \in X} S_x $ \\

We note that $S_x$ plays the role of "fiber" over $x$, while $s = ( U, s, {\cal S} )$ plays
the role of "section" over $U$. \\

Now we define $\pi : {\cal S} \longrightarrow X$ by \\

$ \pi ( S_x ) = \{ x \},~~~ x \in X $ \\

thus it follows that \\

$ {\cal S}_x = \pi^{-1} ( x ) = S_x,~~~ x \in X $ \\

Further, we assume the conditions \\

i)~~~ $ \forall~~ s = ( U, s, {\cal S} ) \in \sigma,~ x \in s ( U ) ~~:~ s ( x ) \in S_x $ \\

ii)~~~ $ \forall~~ x \in X ~:~ S_x \subseteq \bigcup_{s \in \sigma} s ( U ) $ \\

iii)~~~ $ \begin{array}{l}
                   \forall~~ s = ( U, s, {\cal S} ),~ t = ( V, t, {\cal S} ) \in \sigma,~ z \in
                                    s ( U ) \cap t ( V ) ~:~ \\ \\
                   \exists~~ x \in W \subseteq U \cap V,~ W ~~\mbox{open} ~:~ \\ \\
                   ~~~~ z = s ( x ) = t ( x ) \\ \\
                   ~~~~ s|_W = t|_W
           \end{array} $ \\ \\

Clearly, i) is equivalent with \\

$ \forall~~ s ( U, s, {\cal S} ) \in \sigma ~:~ s \in \prod_{x \in U} S_x $ \\

while ii) is equivalent with \\

$ {\cal S} = \bigcup_{s \in \sigma} s ( U ) $ \\

Also, given any sheaf $( {\cal S}, \pi, X )$, if we take \\

$ S_x = {\cal S}_x,~~~ x \in X $ \\

and \\

$ \sigma = \{~ s ~|~ s ~~\mbox{section in}~ ( {\cal S}, \pi, X ) ~\} $ \\

then i), ii) and iii) are obviously satisfied. \\

We also have the converse, namely \\

{\bf Theorem 3.1.} \\

Consider $( {\cal S}, \pi, X )$ as given in (3.1) - (3.3), i) - iii). Then \\

$ {\cal B} = \{~ s ( V ) ~|~  s = ( U, s, {\cal S} ) \in \sigma,~~ V \subseteq U,~ V
                  ~~\mbox{open} ~\} $ \\

is a basis for a topology on ${\cal S}$ for which $( {\cal S}, \pi, X )$ is a sheaf on $X$. \\

Also, for every $s = ( U, s, {\cal S} ) \in \sigma$, the mapping $s : U \longrightarrow
{\cal S}$ is {\it continuous} and {\it open}. \\

Furthermore, ${\cal B}$ gives the {\it finest} topology on ${\cal S}$, for which all the
mappings $s \in\sigma$ are continuous. \\ \\

{\bf 4. EXAMPLES OF GERMS AND SHEAVES} \\

{\bf Sheaf of Germs of Continuous Functions} \\

Let \\

(4.1)~~~ $ {\cal C}_{loc} ( X ) = \sum_{U \in \tau_ X}\, {\cal C} ( U, \mathbb{C} ) $ \\

and for $x \in X$, let \\

(4.2)~~ $ {\cal C}_{loc} ( X )_x =
                            \sum_{x \in U \in \tau_ X}\, {\cal C} ( U, \mathbb{C} ) $ \\

For $x \in X,~ ( U, f, \mathbb{C} ),~ ( V, g, \mathbb{C} ) \in {\cal C}_{loc} ( X )_x$, we
define the equivalence relation on  ${\cal C}_{loc} ( X )_x$ \\

(4.3)~~~ $ f \thicksim_x g ~~~\Longleftrightarrow~~~
                    \left ( \begin{array}{l}
                                     \exists~~ x \in W \in \tau_X,~
                                                      W \subseteq U \cap V ~: \\ \\
                                     ~~~~ f|_W = g|_W
                             \end{array} \right ) $ \\

We denote now \\

(4.4)~~~ $ {\cal C}_{X,x} = {\cal C}_{loc} ( X )_x / \thicksim_x $ \\

while the equivalence class of $( U, f, \mathbb{C} ) \in {\cal C}_{loc} ( X )_x$ is denoted
by \\

(4.5)~~~ $ [ f ]_x = f_x $ \\

and it is called the {\it germ} of $f$ at $x$. \\

We now apply Theorem 3.1., and take \\

(4.6)~~~ $ S_x = {\cal C}_{X,x},~~~ x \in X $ \\

then \\

(4.6)~~~ $ {\cal S} = {\cal C}_X = \sum_{x \in X}\, S_x = \sum_{x \in X}\, {\cal C}_{X,x} $ \\

and define $\pi : {\cal C}_X \longrightarrow X $ by \\

(4.7)~~~ $ \pi( {\cal C}_{X,x} ) = \{ x \},~~~ x \in X $ \\

Further, we take $\sigma$ as follows. For $U \subseteq X$, we have the mapping \\

(4.8)~~~ $ {\cal C}_{loc} ( X ) \ni f ( U, f, \mathbb{C} ) \longmapsto
                                \tilde f = ( U, \tilde f, {\cal C}_X ) \in {\cal C}_X  $ \\

where \\

(4.9)~~~ $ \tilde f ( x ) = [ f ]_x,~~~ x \in X $ \\

and then we take \\

(4.10)~~~ $ \sigma = \{~ \tilde f ~|~ f \in {\cal C}_{loc} ( X ) ~\} $ \\

It is now easy to see that the conditions of Theorem 3.1. are satisfied, thus \\

(4.11)~~~ $ ( {\cal C}_X, \pi, X ) $ ~is a sheaf \\

and thus the set \\

(4.12)~~~ $ {\cal B} = \{~ \tilde f ( V ) ~|~  f \in {\cal C}_{loc} ( X ),~~
                            V \subseteq U,~ V ~~\mbox{open} ~\} $ \\

is a basis for the topology on ${\cal C}_X$ \\

{\bf Remark 4.1.} \\

1) There is a second way to construct the sheaf $( {\cal C}_X, \pi, X )$, namely, starting with
the {\it presheaf} $( {\cal C} ( U, \mathbb{C} ), \rho^U_V )$, presented in section 5. \\

2) The mapping (4.8) is {\it bijective}. \\

3) The above construction can be done for ${\cal C}^l$-smooth functions, with $1 \leq l \leq
\infty$, as well as for analytic functions. \\ \\

{\bf 5. PRESHEAVES} \\

{\bf Definition 5.1. PRESHEAF} \\

$S = ( S ( U ), \rho^U_V )$ is a {\it presheaf} on the topological space $X$, if and
only if \\

(5.1)~~~ $ \tau_X \ni U \longmapsto S ( U ) $ \\

with $S ( U )$ being arbitrary sets, while \\

(5.1a)~~~ $ \rho^U_V : S ( U ) \longrightarrow S ( V ),~~~
                                        U, V \in \tau_X,~~ V \subseteq U $ \\

such that \\

(5.2)~~~ $ \rho^U_U = id_{S ( U )},~~~ u \in \tau_X $ \\

(5.3)~~~ $ \rho^U_W = \rho^V_W \circ \rho^U_V,~~~
                       U, V, W \in \tau_X,~~ W \subseteq V \subseteq U $ \\

{\bf Remark 5.1.} \\

Here we mention two useful alternative interpretations of presheaves. \\

1) Let $\leq$ be the partial order on $\tau_X$ defined by \\

(5.4)~~~ $ U \leq V ~~~\Longleftrightarrow~~~ V \subseteq U $ \\

then $( \tau_X, \leq )$ is right directed and $( S ( U ), \rho^U_V )$ is a right directed
system. This is important when we shall define certain {\it germs}. \\

2) We can see $\tau_X$ as a {\it category}, namely, with the {\it objects} $U \in \tau_X$,
and with the unique {\it morphism} between two objects $V \subseteq U$, namely \\

(5.5)~~~ $ Hom ( V, U ) = \{ \subseteq \} $ \\

Then the presheaf $S = ( S ( U ), \rho^U_V )$ on $X$ is a {\it contravariant functor} from
the category $\tau_X$ to the category $Set$ of all sets and all functions between sets. \\

{\bf Definition 5.2 RESTRICTION OF A PRESHEAF} \\

Given a presheaf $S = ( S ( U ), \rho^U_V )$ on $X$ and an open subset $A \subseteq X$. Then
the {\it restriction} of $S$ to $A$ is the presheaf \\

(5.6)~~~ $ S|_A = ( S ( U ), \rho^U_V ) $ \\

where $V \subseteq U \subseteq A$ are open sets. Thus \\

(5.7)~~~ $ ( S|_A ) ( U ) = S ( U ),~~~ ( \rho|_A )^U_V = \rho^U_V,~~~
                                         V \subseteq U \subseteq A $ open \\

{\bf Remark 5.2.} \\

We have the construction, see further details in section 9 \\ \\

$\begin{array}{||l}
          ~~~SHEAF ~~~------>~~ PRESHEAF ~ OF ~ SECTIONS \\ \\
          ~~~{\cal S} = ( {\cal S}, \pi, X ) ~~------>~~
                    \Gamma ( {\cal S} ) = ( \Gamma ( U, {\cal S} ), \rho^U_V = |_V )
 \end{array} $ \\ \\ \\

 {\bf 6. MORPHISMS OF PRESHEAVES} \\

 Given two presheaves $S = ( S ( U ), \rho^U_V )$ and $E = ( E ( U ), \lambda^U_V )$ on $X$,
 by {\it morphism} from $S$ to $E$, that is \\

(6.1)~~~ $ \phi : S \longrightarrow E $ \\

one means a family of mappings \\

(6.2)~~~ $ \phi = ( \phi_U )_{U \in \tau_X} $ \\

with the following two properties \\

(6.3)~~~ $ \phi_U : S ( U ) \longrightarrow E ( U ),~~~ U \in \tau_X $ \\

and with the commutative diagram \\

\begin{math}
\setlength{\unitlength}{0.2cm}
\thicklines
\begin{picture}(60,20)

\put(10,16){$S ( U )$}
\put(28,18){$\phi_U$}
\put(16,16.5){\vector(1,0){28.5}}
\put(46,16){$E ( U )$}
\put(0,9.5){$(6.4)$}
\put(12,14.5){\vector(0,-1){9}}
\put(8.5,9.5){$\rho^U_V$}
\put(48,14.5){\vector(0,-1){9}}
\put(49,9.5){$\lambda^U_V$}
\put(10,2){$S ( V )$}
\put(16,2.6){\vector(1,0){28.5}}
\put(46.5,2){$E ( V )$}
\put(27.5,-0.5){$\phi_V$}

\end{picture}
\end{math} \\

for open subsets $ V \subseteq U \subseteq X$. \\

Now we can define the {\it category} of presheaves on $X$, namely \\

(6.5)~~~ $ {\cal PS}h_X $ \\

as having {\it objects} all the presheaves on $X$, and with the {\it morphisms} being all
the presheaf morphism between any two presheaves on $X$. \\

{\bf Definition 6.1.} \\

Given a morphism $\phi = ( \phi_U )_{U \in \tau_X} : S \longrightarrow E$ in the category
${\cal PS}h_X$, then $\phi$ is an {\it injection}, if and only if each $\phi_U : S ( U )
\longrightarrow E ( U )$ is injective. Similarly,  $\phi$ is a {\it surjection}, if and only
if each $\phi_U : S ( U ) \longrightarrow E ( U )$ is surjective. \\
Also,  $\phi$ is an {\it isomorphism}, if and only if each $\phi_U : S ( U ) \longrightarrow
E ( U )$ is a bijection. \\ \\

{\bf 7. FIBERS OF PRESHEAVES, SHEAFIFICATION, \\
        \hspace*{0.4cm} SHEAF GENERATED BY A PRESHEAF} \\ \\

$\begin{array}{||l}
          ~~~PRESHEAF ~~------>~~ SHEAF \\ \\
          ~~~ S ~~------>~~ \mathbb{S} ( S )
 \end{array} $ \\ \\

Given a presheaf $S = ( S ( U ), \rho^U_V )$ on $X$, with the help of Theorem 3.1., we shall
associate with it a sheaf $\mathbb{S} ( S )$ on $X$, called the {\it sheafification} of
$S$. \\

Let $x \in X$ and ${\cal V}( x )$ be the set of open neighbourhood of $x$. Let ${\cal B}( x )$
be any basis of open neighbourhoods of $x$. Then $( S ( U ), \rho^U_V )$, with $V, U \in
{\cal V}( x )$, as well as $( S ( U ), \rho^U_V )$, with $V, U \in {\cal B}( x )$ are right
directed systems. Furthermore, ${\cal B}( x )$ is cofinal in ${\cal V}( x )$. Therefore, we
can define the {\it fiber} of the presheaf $S$ over $x$, by \\

(7.1)~~~ $ S_x = \underrightarrow{\lim}_{\,U \in {\cal V}( x )\,} S ( U ) \simeq
                      \underrightarrow{ \lim }_{\,U \in {\cal B}( x )\,} S ( U ) $ \\

where the second relation is a bijection of sets. \\

Now according to Theorem 3.1., we define \\

(7.2)~~~ $ {\cal S} = \sum_{x \in X} S_x $ \\

(7.3)~~~ $ \pi : {\cal S} \longrightarrow X $ \\

by \\

(7.4)~~~ $ \pi ( S_ x ) = \{ x \},~~~ x \in X$ \\

Thus \\

(7.5)~~~ $ {\cal S}_x = S_x = \pi^{-1} ( x ),~~~ x \in X $ \\

It only remains to define the topology on ${\cal S}$ through the family of mappings \\

(7.6)~~~ $ \sigma = \{~ \tilde s = ( U, \tilde s, {\cal S} ) ~\},~~~
                               U \in \tau_X,~~ \tilde s : U \longrightarrow {\cal S} $ \\

which is done next in several steps. \\

First, for open $U \subseteq X$, we define the mapping $\rho_U$ by, see further details in
(7.15) \\

(7.7)~~~ $ S ( U ) \ni s \longmapsto \rho_U ( s ) = \tilde s $ \\

where \\

(7.8)~~~ $ \tilde s : U \longrightarrow {\cal S} $ \\

with \\

(7.9)~~~ $ U \ni x \longmapsto \rho^x_U ( s ) \in S_x = {\cal S}_x \subseteq {\cal S} $ \\

while, for $x \in X$, we have \\

(7.10)~~~ $ \rho^x_U : S ( U ) \longrightarrow S_x =
                      \underrightarrow{\lim}_{\,U \in {\cal V}( x )\,} S ( U ) $ \\

which is the {\it canonical} mapping corresponding to the direct limit in the right hand
side, see (A.8) in Appendix. For convenience, we shall also use the notation \\

(7.11)~~~ $ \tilde s ( x ) = ( \rho_U ( s ) ) ( x ) = \rho^x_U ( s ) = [ s ]_x,~~~
                                         U \in \tau_X,~~ s \in S ( U ), ~~ x \in U $ \\

And now, we have obtained (7.6) in the form \\

(7.12)~~~ $ \sigma = \{~ \tilde s = ( U, \tilde s, {\cal S} ) ~~|~~
                                           U \in \tau_X,~ s : S ( U ) ~\} $ \\

Applying Theorem 3.3., it follows that \\

(7.13)~~~ $ \mathbb{S} ( S ) = \mathbb{S} ( ( S ( U ), \rho^U_V ) ) =
                                             {\cal S} = ( {\cal S}, \pi, X ) $ \\

is a sheaf on $X$, called the {\it sheafification} of $S = ( ( S ( U ), \rho^U_V )$. \\

Let us end by noting the following useful relation, see (A.12) in Appendix \\

(7.14)~~~ $ \rho^x_U = \rho^x_V \circ \rho^U_V,~~~
                           x \in V \subseteq U \subseteq X,~~ U, V $ open \\

Furthermore, now that ${\cal S} = ( {\cal S}, \pi, X )$ was shown to be a sheaf, we can
specify the {\it range} of the mappings $\rho_U$ in (7.7), namely \\

(7.15)~~~ $ \rho_U : S ( U ) \longrightarrow {\cal S} ( U ) =
                                      \Gamma ( U, {\cal S} ),~~~ U \in \tau_X $ \\

thus these mappings associate individual sets $S ( U )$ of the presheaf $S$ with sections
${\cal S} ( U )$ of the sheaf ${\cal S} = \mathbb{S} ( S )$. \\ \\

{\bf 8. THE SHEAFIFICATION FUNCTOR} \\

It is easy to see that \\

(8.1)~~~ $ \mathbb{S} : {\cal PS}h_X \longrightarrow {\cal S}h_X $ \\

is a {\it covariant functor}. \\

In fact, as seen in Theorem 9.1. below, when restricted to the subcategory of {\it complete
presheaves} of ${\cal PS}h_X$, the functor $\mathbb{S}$ has as {\it right adjoint} the functor
$\Gamma$ of {\it sections}. \\ \\

{\bf 9. PRESHEAF OF SECTIONS OF A SHEAF. \\
        \hspace*{0.4cm} THE SECTION FUNCTOR} \\ \\

$\begin{array}{||l}
          ~~~ SHEAF ~~------>~~ PRESHEAF ~ OF ~ SECTIONS ~~------>~~ \\ \\
          ~~~~~~~~~~~~------>~~ SHEAFIFICATION ~=~ ORIGINAL ~ SHEAF\\ \\
          ~~~ {\cal S} ~~------>~~ \Gamma ( {\cal S} ) ~~------>~~ \mathbb{S} ( \Gamma ( {\cal S} ) ) = {\cal S}
 \end{array} $ \\ \\

Here we further develop the construction in Remark 5.2. \\

Given any {\it sheaf} \\

(9.1)~~~ $ {\cal S} = ( {\cal S}, \pi, X ) $ \\

on a topological space $X$, we shall associate with it its {\it presheaf} of sections. This
construction can be done in the general setting of functors, and in this case, of the {\it
section functor} \\

(9.2)~~~ $ \Gamma ( . , {\cal S} ) = \Gamma_{\cal S} : \tau_X \longrightarrow Set $ \\

namely, regarding the mappings between {\it objects} in the two categories we have \\

(9.3)~~~ $ \Gamma ( . , {\cal S} ) ( U ) = \Gamma_{\cal S} ( U ) =
                              \Gamma ( U, {\cal S} ) = {\cal S} ( U ),~~~ U \in \tau_X $ \\

while the mappings between {\it morphisms} in the two categories are \\

(9.4)~~~ $ \sigma^U_V : \Gamma ( U, {\cal S} ) \longrightarrow
                             \Gamma ( V, {\cal S} ),~~~ U, V \in \tau_X,~~ V \subseteq U $ \\

where \\

(9.5)~~~ $ \sigma^U_V ( s ) = s|_V,~~~ s \in \Gamma ( U, {\cal S} ) $ \\

with $s|_V$ being the usual restriction to $V$ of the function $s : U
\longrightarrow \mathbb{C}$. \\

In this way, the {\it sheaf} ${\cal S}$ on $X$ in (9.1) is associated with the {\it presheaf}
on $X$ \\

(9.6)~~~ $ \Gamma ( {\cal S} ) = \Gamma_{\cal S} = ( \Gamma ( U, {\cal S} ), \sigma^U_V ) $ \\

We recall that $( \Gamma ( U, {\cal S} ), \sigma^U_V )$ is a right directed system, thus \\

(9.7)~~~ $ \Gamma ( {\cal S} )_x = ( \Gamma_{\cal S} )_x =
           \underrightarrow{\lim}_{\,x \in U \in \tau_X\,} \Gamma ( U, {\cal S} ),~~~ x \in X $ \\

As a consequence, we obtain \\

(9.8)~~~ $ {\cal S}_x \backsimeq ( \Gamma_{\cal S} )_x,~~~ x \in X $ \\

Also, let $S = ( S ( U ), \rho^U_V )$ be a {\it presheaf} on the topological space $X$, let
${\cal S} = \mathbb{S} ( S )$, and let $\Gamma_{\cal S} = ( \Gamma ( U, {\cal S} ),
\sigma^U_V )$, then \\

(9.9)~~~ $ S_x \backsimeq ( \Gamma_{\cal S} )_x,~~~ x \in X $ \\

The main result regarding the functors of sheafification and sections is given in \\

{\bf Theorem 9.1.} $\Gamma$ {\bf RIGHT ADJOINT TO} $\mathbb{S}$ \\

Let ${\cal S} = ( {\cal S}, \pi, X )$ be a sheaf on a topological space $X$, then \\

(9.10)~~~ $ {\cal S} = \mathbb{S} ( \Gamma_{\cal S} ) $ \\

that is \\

(9.11)~~~ $ {\cal S} \stackrel{\Gamma}\longmapsto \Gamma_{\cal S}
                          \stackrel{\mathbb{S}}\longmapsto {\cal S} $

\hfill $\Box$ \\

For presheaves $S = ( S ( U ), \rho^U_V )$ on $X$, one also has, for every open $A \subseteq
X$, the restriction property of sheafification \\

(9.12)~~~ $ \mathbb{S} ( S|_A ) = ( \mathbb{S} ( S) )|_A $ \\ \\

{\bf 11. COMPLETE PRESHEAVES : \\
         \hspace*{0.7cm} THE ORIGINAL 1946 CONCEPT OF JEAN LERAY} \\

{\bf Definition 11.1. COMPLETE PRESHEAF} \\

A presheaf $S = ( S ( U ), \rho^U_V )$ on a topological space $X$ is called  {\it complete},
if and only if the following two conditions hold \\

(S1)~~~ $ \begin{array}{||l}
             ~~\forall~~ U \in \tau_X,~ ( U_\alpha )_{\alpha \in I} ~~\mbox{open cover of}~
                           U,~~ s, t \in S ( U ) ~: \\ \\
             ~~~~~~ (~ \rho^U_{U_\alpha} ( s ) = \rho^U_{U_\alpha} ( t ),~~ \alpha \in I ~)
                       ~~\Longrightarrow~~  s = t
           \end{array} $ \\ \\ \\

(S2)~~~ $ \begin{array}{||l}
             ~~\forall~~ U \in \tau_X,~ ( U_\alpha )_{\alpha \in I} ~~\mbox{open cover of}~
                           U,~ ( s_\alpha )_{\alpha \in I} \in
                                     \prod_{\alpha \in I} S ( U_\alpha ) ~: \\ \\
             ~~~~ \left ( \begin{array}{l}
                               \forall~~ \alpha, \beta \in I ~: \\ \\
                               ~~ U_{\alpha,\beta} = U_\alpha \cap U_\beta \neq \phi
                                    \Longrightarrow s_\alpha|_{U_{\alpha,\beta}} =
                                           s_\beta|_{U_{\alpha,\beta}}
                             \end{array} \right ) \Longrightarrow
                                 \left ( \begin{array}{l}
                                                \exists~~ s \in S ( U ) ~: \\ \\
                                                ~~ s|_{U_\alpha} = s_\alpha,~ \alpha \in I
                                            \end{array} \right )
          \end{array} $ \\ \\

Here we note that (S1) implies the {\it uniqueness} of $s$ in (S2). \\

{\bf Lemma 11.1. PRESHEAVES OF SECTIONS \\
                 \hspace*{2.75cm} ARE COMPLETE} \\

Let ${\cal S}$ be a sheaf on the topological space $X$, then its presheaf of sections
$\Gamma_{\cal S}$ is complete, that is \\

(11.1)~~~ $ {\cal S} \stackrel{\Gamma}\longmapsto \Gamma_{\cal S}
                                           ~~\mbox{complete presheaf} $ \\

{\bf Proposition 11.1. CHARACTERIZATION OF \\
                       \hspace*{3.65cm} COMPLETE PRESHEAVES} \\

Let $S = ( S ( U ), \rho^U_V )$ be a presheaf on the topological space $X$, and let ${\cal S}
= \mathbb{S} ( S )$. Then we have the equivalent properties, see (7.15) \\

(11.2)~~~ $ S ~~\mbox{is a complete presheaf} $ \\

(11.3)~~~ $ \rho_U : S ( U ) \stackrel{bijective}\longrightarrow
                  \Gamma ( U, {\cal S} ) = {\cal S} ( U ),~~~ U \in \tau_X $

\hfill $\Box$ \\

In particular \\

(11.4)~~~ $ S ~~\mbox{satisfies (S1)} ~~~\Longleftrightarrow~~~
                                \mbox{the mappings (11.3) are injective} $ \\

Also \\

(11.5)~~~ $ S ~~\mbox{complete} ~~~ \Longrightarrow~~~
                             (~ S|_A ~~\mbox{complete},~~ A \in \tau_X ~) $ \\

In particular, for the presheaf of sections of $\Gamma( {\cal S} )$ of a sheaf ${\cal S}$, we
have \\

(11.6)~~~ $ \Gamma ( {\cal S}|_A ) = ( \Gamma ( {\cal S} ) )|_A,~~ A \in \tau_X $ \\

A stronger version of the characterization of complete presheaves in Proposition 11.1 is given
in \\

{\bf Theorem 11.2. CHARACTERIZATION OF \\
                   \hspace*{3.1cm} COMPLETE PRESHEAVES} \\

Let $S$ be a presheaf on the topological space $X$, then we have the equivalent properties \\

(11.7)~~~ $ S ~~\mbox{is a complete presheaf} $ \\

(11.8)~~~ $ S = \Gamma ( \mathbb{S} ( S ) ) $ \\

(11.9)~~~ $ \exists~~ {\cal E} ~~\mbox{sheaf on}~ X ~:~ S = \Gamma ( {\cal E} ) $

\hfill $\Box$ \\

Furthermore, $\Gamma ( \mathbb{S} ( S ) )$ turns out to be the {\it smallest} complete
presheaf containing $S$, as shown in \\

{\bf Theorem 11.3. SMALLEST COMPLETE PRESHEAF} \\

Let $S = ( S ( U ), \rho^U_V )$ be a presheaf on the topological space $X$, ${\cal S} =
\mathbb{S} ( S )$ its sheafification, and $\Gamma ( {\cal S} ) = ( \Gamma ( U, {\cal S} ),
\sigma^U_V )$ the presheaf of sections of ${\cal S}$. Then \\

(11.10)~~~ $ \begin{array}{l}
                     \forall~~ E ~~\mbox{complete presheaf on}~ X,~~ \phi : S \longrightarrow
                                         E ~~\mbox{presheaf morphism} ~: \\ \\
                     \exists~!~~ \psi : \Gamma ( {\cal S} ) \longrightarrow E
                         ~~\mbox{presheaf morphims with the commutative diagram} ~:
            \end{array} $ \\ \\

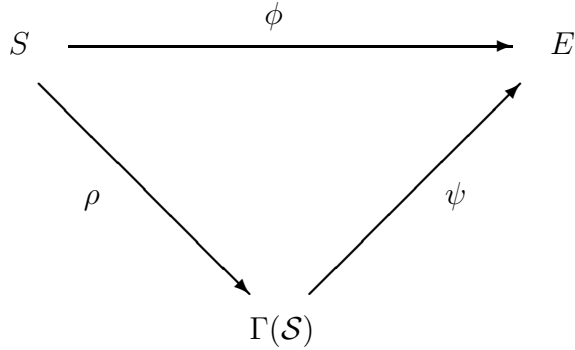
\begin{math}
\setlength{\unitlength}{0.2cm}
\thicklines
\begin{picture}(60,20)

\put(16,16){$S$}
\put(33,18){$\phi$}
\put(20,16.5){\vector(1,0){29.5}}
\put(52,16){$E$}
\put(21,6){$\rho$}
\put(18,14){\vector(1,-1){14}}
\put(32,-3){$\Gamma ( {\cal S} )$}
\put(45,6){$\psi$}
\put(36,0){\vector(1,1){14}}

\end{picture}
\end{math} \\ \\

where the presheaf morphism \\

(11.11)~~~ $ \rho = ( \rho_U )_{U \in \tau_X} : S \longrightarrow \Gamma ( {\cal S} ) $ \\

is given by (11.3). \\

{\bf Remark 11.1. ON COMPLETENESS OF PRESHEAVES} \\

A presheaf $S = ( S ( U ), \rho^U_V )$ on a topological space $X$ is called {\it functional},
if and only if \\

1) the sections $s \in S ( U )$ are usual functions defined on $U$ \\

2) the sheaf morphisms $\rho^U_V : S ( U ) \longrightarrow S ( V )$ are usual restrictions
of \\
\hspace*{0.4cm} functions. \\

Obviously \\

(11.12)~~~ every functional presheaf satisfies (S1) \\

A functional presheaf $S = ( S ( U ), \rho^U_V )$ on a topological space $X$ is called {\it
local}, if and only if the sections $s \in S ( U )$ are defined by local conditions, such as
for instance, continuity, various levels of smoothness, analyticity, local integrability,
etc. \\

As an example of {\it nonlocal} condition is {\it boundedness}. \\

Now, it is easy to see that \\

(11.13)~~~ every local presheaf satisfies both (S1) and (S2), thus it is \\
           \hspace*{1.7cm} complete \\

Such were what, back in 1946, Jean Leray first considered to be "sheaves". \\

A {\it counterexample} : let $X$ be any infinite set with the discrete topology, and let
$S = ( S ( U ), \rho^U_V )$ be the functional presheaf on $X$ with \\

(11.14)~~~ $S ( U ) = \{~ s : U \longrightarrow \mathbb{R} ~~|~~ s
                                        ~~\mbox{bounded on}~ U ~\} $ \\

which obviously satisfies (S1). However \\

(11.15)~~ $ S $ is {\it not} local \\

Indeed, any $ s \in S ( U )$, with $U$ an infinite subset of $X$, is necessarily bounded at
every $x \in U$, since it has a finite value $s ( x ) \in \mathbb{R}$. However, it obviously
need not be bounded on the whole of $U$. \\

And it is easy to see that \\

(11.16)~~~ the functional presheaf $ S $ in (11.14) is {\it not} complete \\

since it fails to satisfy (S2). \\

We denote by \\

(11.17)~~~ $ {\cal C}o{\cal PS}h_X $ \\

the category of {\it complete presheaves} on the topological space $X$. \\ \\

{\bf 13. THE SECTION FUNCTOR $\Gamma$ \\
         \hspace*{0.65cm} AND ITS LEFT ADJOINT, \\
         \hspace*{0.65cm} THE SHEAFIFICATION FUNCTOR $\mathbb{S}$} \\ \\

$\begin{array}{||l}
           ~~~ {\cal S}h_X ~ \stackrel{\Gamma}\simeq ~ {\cal C}o{\cal PS}h_X ~~ ARE ~
                              EQUIVALENT ~ CATEGORIES
 \end{array} $ \\ \\

{\bf Theorem 13.1. EQUIVALENCE OF SHEAVES AND \\
                   \hspace*{3.1cm} COMPLETE PRESHEAVES} \\

The mapping \\

(13.1)~~~ $ {\cal S}h_X \ni {\cal S} \stackrel{\Gamma}\longrightarrow
                        \Gamma ( {\cal S} ) \in {\cal C}o{\cal PS}h_X $ \\

is a {\it covariant} functor which is an {\it equivalence} of the respective two
categories. \\

Specifically \\

(13.2)~~~ $ \mathbb{S}\, \Gamma = id_{\,{\cal S}h_X},~~~~~~
                        \Gamma\, \mathbb{S} = id_{\,{\cal C}o{\cal PS}h_X} $ \\ \\

{\bf 14. CHANGE OF BASE SPACE} \\

Let $X, Y$ be two topological spaces together with a continuous mapping \\

(14.1)~~~ $ f : X \longrightarrow Y $ \\

{\bf 14.1. PUSH-OUT FUNCTORS ON PRESHEAVES \\
           \hspace*{1cm} AND SHEAVES} \\

First we define the {\it covariant} functor, called the {\it push-out} of $f$, namely \\

(14.2)~~~$ \tilde f_* : {\cal PS}h_X \longrightarrow {\cal PS}h_Y $ \\

as follows. Given $S = ( S ( U ), \rho^U_V ) \in {\cal PS}h_X$ and $V \in \tau_Y$, we
define \\

(14.3)~~~ $ ( \tilde f_* ( S ) ) ( V ) = S ( f^{-1} ( V ) ) $ \\

Then for every presheaf morphism $\phi : S \longrightarrow S\,'$ in ${\cal PS}h_X $, one can
define the corresponding presheaf morphism in  ${\cal PS}h_Y $, namely \\

(14.4)~~~ $ \tilde f_* ( \phi ) : \tilde f_* ( S ) \longrightarrow \tilde f_* ( S\,' ) $ \\

given by, see (6.1) - (6.4) \\

(14.5)~~ $ ( \tilde f_* \phi ) )_V = \phi_{f^{-1} ( V )},~~~ V \in \tau_Y $ \\

An important property of the push-out functor is \\

(14.6)~~~ $ \tilde f_* : {\cal C}o{\cal PS}h_X \longrightarrow {\cal C}o{\cal PS}h_Y $ \\

that is, complete presheaves $S$ on $X$ are taken into complete presheaves $\tilde f_* ( S )$
on $Y$. \\

Now, as a main interest, we define the corresponding push-out functors for sheaves, namely \\

(14.7)~~~ ~$ f_* : {\cal S}h_X \longrightarrow {\cal S}h_Y $ \\

and do so simply by the following composition of functors \\

(14.8)~~~ $ f_* = \mathbb{S}_Y \circ \tilde f_* \circ \Gamma_X $ \\

which is a correct definition, since as seen above, we have \\

(14.9)~~~ $ {\cal S}h_X \stackrel{\Gamma_X}\longrightarrow {\cal PS}h_X
                          \stackrel{\tilde f_*}\longrightarrow {\cal PS}h_Y
                            \stackrel{\mathbb{S}_Y}\longrightarrow {\cal S}h_Y $ \\

{\bf 14.2. PULL-BACK FUNCTORS ON SHEAVES} \\

Now we define the {\it covariant} functor, called the {\it pull-back} of $f$, namely \\

(14.10)~~~$ f^* : {\cal S}h_Y \longrightarrow {\cal S}h_X $ \\

as follows. Given ${\cal E} = ( {\cal E}, \pi, Y ) \in {\cal S}h_Y$, then we define \\

(14.11)~~~ $ f^* ( {\cal E} ) = {\cal S} = ( {\cal S}, \lambda, X ) $ \\

where \\

(14.12)~~~ $ {\cal S} = \{~ ( x, z ) \in X \times {\cal E} ~~|~~ f ( x ) = \pi ( z ) ~\} $ \\

while \\

(14.13)~~~ $ \lambda : {\cal S} \longrightarrow X $ \\

is given by \\

(14.14)~~~ $ \lambda ( x, z ) = z,~~~ ( x, z ) \in {\cal S} $ \\

An alternative definition of the pull-back functor is as follows. Let ${\cal E} =
( {\cal E}, \pi, Y ) \in {\cal S}h_Y$. We define the {\it complete presheaf} on $X$, given
by \\

(14.15)~~~ $ {\cal E}_f = ( {\cal E}_f ( U ), \nu^U_V ) $ \\

where \\

(14.16)~~~ $ {\cal E}_f ( U ) = \{~ t : U \longrightarrow {\cal E} ~~|~~
                                     t ~~\mbox{continuous}~, \pi \circ t = f ~\} $ \\

while \\

(14.17)~~~ $ \nu^U_V : {\cal E}_f ( U ) \longrightarrow {\cal E}_f ( V ) $ \\

is the usual restriction of functions. And now we obtain \\

(14.18)~~~  $ f^* ( {\cal E} ) = \mathbb{S} ( {\cal E}_f ) $ \\

as an alternative definition to (14.11). \\ \\

\newpage

{\bf \Large COMMENTS ON : PRESHEAVES, SHEAVES, \\ \\ BUNDLES, MANIFOLDS \\ \\ and SINGULARITIES} \\ \\

It is useful to compare the concepts of presheaf, sheaf, bundle and
manifold, and see their differences not only in their different
formal mathematical structures, but also in the very
purposes for which they were introduced in the first place. \\ \\

{\bf 1. MANIFOLDS} \\

The oldest among the mentioned concepts, namely, that of manifold,
was introduced by Riemann, and it has had an exceptionally important
role not only in differential geometry, but also in
various theories of physics, among them general relativity. \\
The aim of the concept of manifold is to deal with geometric spaces
which are no longer "flat" like the Euclidean ones, and instead have
"curvature". The inevitable effect of such a curvature is that such
spaces can no longer be modelled globally by Euclidean spaces, and
instead, all one can hope for as the next best, is to model them
locally with such flat Euclidean spaces. Consequently, in a modern
formulation, much different from that of Riemann,
we have \\

{\bf Definition 1.1. MANIFOLDS} \\

A topological space $X$ is an $n$-{\it dimensional smooth manifold},
or in short, an $n$-{\it manifold}, if and only if for every $x \in
X$, there exists a neigbourhood $U \subseteq X$ of $x$ and a
homeomorphism $\psi_U : U \longrightarrow \mathbb{R}^n$, with the
following property
of compatibility on overlaps of such neighbourhoods : \\

Given two homeomorphisms $\psi_U : U \longrightarrow \mathbb{R}^n$
and $\psi_V : V
\longrightarrow \mathbb{R}^n$, with $U \cap V \neq \phi$, then the mapping \\

(1.1)~~~ $ \psi_V \circ \psi_U^{-1} :
                 \psi_U ( U \cap V ) \longrightarrow \psi_V ( U \cap V ) $ \\

is a usual smooth mapping of open subsets of $\mathbb{R}^n$. \\

Clearly, we have here the commutative diagram of mappings \\

\begin{math}
\setlength{\unitlength}{0.2cm} \thicklines
\begin{picture}(60,20)

\put(11,16){$U \cap V$} \put(6.3,16){$X \supseteq$}
\put(30,18){$id_{U \cap V}$} \put(18,16.5){\vector(1,0){28.5}}
\put(49,16){$U \cap V$} \put(55,16){$\subseteq X$}
\put(0,9.5){$(1.2)$} \put(14,14.5){\vector(0,-1){9}}
\put(10.5,9.5){$\psi_U$} \put(51.5,14.5){\vector(0,-1){9}}
\put(52.5,9.5){$\psi_V$} \put(8.5,2){$\psi_U ( U \cap V )$}
\put(3,2){$\mathbb{R}^n \supseteq$} \put(20,2.6){\vector(1,0){24}}
\put(46,2){$\psi_V ( U \cap V )$} \put(29.5,-0.5){$\psi_V \circ
\psi_U^{-1}$} \put(56,2){$\,\subseteq\, \mathbb{R}^n$}

\end{picture}
\end{math} \\

Any pair $( U, \psi_U )$ is called a {\it chart} of the manifold
$X$, while a collection
$( U_\alpha, \psi_{U_\alpha} )_{\alpha \in I}$ of charts for which \\

(1.3)~~~ $ X = \bigcup_{\alpha \in I} U_\alpha $ \\

is called an {\it atlas} of the manifold $X$. \\

Obviously, a manifold can in general have many atlases, thus, even a
larger amount of
charts. \\

In case a manifold has an atlas with one single chart, then it is
trivial, since the space $X$ is diffeomorphic with $\mathbb{R}^n$,
see Example 1.1. below. Thus manifold theory is
nontrivial only in the case of spaces which do not have an atlas with one single chart. \\

It follows that the price we have to pay for dealing with
"curvature" is to consider two topological spaces, namely, $X$ and
$\mathbb{R}^n$, and on top of that, also atlases,
therefore, quite likely, a multiplicity of charts. \\

{\bf Example 1.1.} \\

A trivial example of $n$-manifold is $\mathbb{R}^n$. Indeed, in this
case we can have an atlas formed from one single chart, namely $\psi
= id_{\mathbb{R}^n} : U = \mathbb{R}^n
\longrightarrow \mathbb{R}^n$. \\

{\bf Example 1.2.} \\

A simple, yet nontrivial example of $1$-manifold is given by the
unit circle in
$\mathbb{R}^2$, namely \\

(1.4) $~~~ S^1 = \{~ x = ( x_1,x_2 ) \in \mathbb{R}^2 ~|~
                      x_1^2 + x_2^2 = 1 ~\} \subseteq \mathbb{R}^2 $ \\

Indeed, in this case every atlas must contain at least two charts.
For instance, one such
chart $( \psi_U, U )$ is given by \\

(1.5) $~~~ U = \{~ ( x, y ) \in \mathbb{R}^2 ~|~ x^2 + y^2 = 1,~~ y > - 1 ~\} \subseteq S^1 $ \\

and a corresponding $\psi_U : U \longrightarrow \mathbb{R}$, while
the second chart
$( \psi_V, V )$ may have \\

(1.6) $~~~ V = \{~ ( x, y ) \in \mathbb{R}^2 ~|~ x^2 + y^2 = 1,~~ y < 1 ~\} \subseteq S^1 $ \\

with the appropriate  $\psi_V : V \longrightarrow \mathbb{R}$. \\

{\bf Example 1.3.} \\

In this regard, the situation does not get more complicated in the
case of the $n$-dimensional
sphere \\

(1.7)~~~ $ S^n = \{~ x = ( x_1, \ldots , x_{n + 1} ) \in
\mathbb{R}^{n + 1} ~
                    x_1^2 + \ldots + x_{n + 1}^2 = 1 ~\} \subseteq \mathbb{R}^{n + 1}$ \\

when again there is an atlas formed by only two charts, namely \\

(1.8) $~~~ U = \left \{ x = ( x_1, \ldots , x_{n + 1} ) \in
\mathbb{R}^{n + 1} ~
                         \begin{array}{|l}
               x_1^2 + \ldots + x_{n + 1}^2 = 1 \\ \\ x_{n + 1} > - 1
                          \end{array} \right \} \subseteq \\ \\
                          \hspace*{11cm} \subseteq S^n $ \\

and a corresponding $\psi_U : U \longrightarrow \mathbb{R}^n$, while
the second chart
$( \psi_V, V )$ may have \\

(1.9)  $~~~ V = \left \{ x = ( x_1, \ldots , x_{n + 1} ) \in
\mathbb{R}^{n + 1} ~
                         \begin{array}{|l}
               x_1^2 + \ldots + x_{n + 1}^2 = 1 \\ \\ x_{n + 1} < 1
                          \end{array} \right \} \subseteq \\ \\
                          \hspace*{11cm} \subseteq S^n $ \\

with the appropriate  $\psi_V : V \longrightarrow \mathbb{R}^n$. \\

{\bf Example 1.4.} \\

Another nontrivial example of manifold, one that is obviously not
globally "flat", is the Moebius band. Indeed, it is a $2$-manifold,
since every point on it has an open neighbourhood diffeomorphic with
$\mathbb{R}$.

\hfill $\Box$ \\

The fact to note is that, unlike with sheaves and bundles, each of
which are defined by two topological spaces and a surjective
continuous mapping between them with certain specific properties, in
the case of manifolds one does not have a similar simplicity, as one
has in addition to the two spaces $X$ and $\mathbb{R}^n$, also a
considerable amount, often an infinity in fact, of charts in an
atlas. And in general, there is no canonical way to choose
the charts, and therefore, the atlases. \\

In this way, the curvature of the space $X$ can be dealt with
locally in a suitable
neighbourhood of each point $x \in X$ with the help of a flat Euclidean space. \\
However, globally, that is for the whole of $X$, dealing with the
curvature leads to the
complexity of many charts, and also atlases. \\ \\

{\bf 2. SHEAVES} \\

{\bf Definition 2.1. SHEAVES} \\

We call $( {\cal S}, \pi, X )$ a {\it sheaf on the base space} $X$,
if and only if $\pi : {\cal S} \longrightarrow X$ is a continuous
surjective mapping which in addition is also a
local homeomorphism. \\

As seen from the examples next, the concept of sheaf is very much
different from that of
manifold, even if both involve essentially local homeomorphisms. \\

{\bf Example 2.1.} \\

A trivial example of sheaf is the {\it constant sheaf} \\

(2.1)~~~ $ {\cal S} = ( {\cal S} = X \times M, pr_X, X) $ \\

where $X$ is a topological space, while $M$ is any set considered
with the discrete topology.
In this case, given any open subset $U \subseteq X$, it is obvious that \\

(2.2)~~~ $ A = U \times \{ m \} \subseteq {\cal S} = X \times M $ \\

is an open subset of ${\cal S} = X \times M$, for every $m \in M$,
and furthermore, $\pi$
restricted to it, that is \\

(2.3)~~~ $ \pi|_A : U \times \{ m \} \ni ( x, m ) \longmapsto x \in U $ \\

is a homeomorphism. \\

{\bf Example 2.2.} \\

A simple nontrivial example of sheaf is obtained as follows. Let
${\cal S} = \mathbb{R}$ and
$X = S^1$, with the mapping $\pi : {\cal S} \longrightarrow X$ defined by \\

(2.4)~~~ $ {\cal S} = \mathbb{R} \ni \alpha \longmapsto
                      x = ( \cos \alpha, \sin \alpha ) \in X = S^1 $ \\

Then for every open and strict subset $U \subsetneqq X = S^1$, there
exists a countable
infinity of open strict subsets $A \subsetneqq {\cal S} = \mathbb{R}$ such that \\

(2.5)~~~ $ \pi|_A : A \longrightarrow U $ \\

is a homeomorphism. In fact, for every open interval $A = ( a, b )
\subsetneqq {\cal S} =
\mathbb{R}$, with $b - a < 2 \pi$, we have the homeomorphism \\

(2.6)~~~ $ \pi|_A : A \longrightarrow \pi ( A ) $

\hfill $\Box$ \\

Similar with bundles, see section 3 below, an important concept in a
sheaf $( {\cal S}, \pi, X )$ is that of {\it fiber}, or {\it stalk}.
Namely, for $x \in X$, the fiber over $x$ is
given by \\

(2.7)~~~ $ {\cal S}_x = \pi^{-1} ( x ) $ \\

The essential difference when compared with bundles is that, in the
case of sheaves, we always
have \\

(2.8)~~~ $ {\cal S}_x ~~\mbox{is a discrete subset of}~ {\cal S} $ \\

for each given $x \in X$. For instance, Example 2.2., when $X =
S^1$, if we take any $x \in
X$, then \\

(2.9)~~~ $ {\cal S}_x = \{ \alpha + 2 k \pi ~|~ k \in \mathbb{Z} \} $ \\

where $\alpha \in {\cal S} = \mathbb{R}$ is such that $x = ( \cos
\alpha, \sin \alpha )$. And
clearly, ${\cal S}_x $ in (2.9) is a discrete subset of ${\cal S} = \mathbb{R}$. \\

In the case of the constant sheaf, we have for $x \in X$ \\

(2.10)~~~ $ {\cal S}_x = \pi^{-1} ( x ) = \{ x \} \times M
                                \subseteq {\cal S} = X \times M $ \\

which again is a discrete subset of ${\cal S} = X \times M$. \\ \\

{\bf 3. BUNDLES} \\

{\bf 3. Definition 3.1. BUNDLES} \\

We call $( E, \pi, X )$ a {\it bundle over the base space} $B$, if
and only if $E$ and $X$ are topological spaces and $\pi : E
\longrightarrow X$ is a continuous surjective mapping. In
this case $E$ is called the {\it bundle space} or {\it total space}. \\

For $x \in X$, the {\it fiber} or {\it stalk over} $x$ is given by \\

(3.1)~~~ $ E_x = \pi^{-1} ( x ) $ \\

A useful particular case of bundles are the {\it fiber bundles}
given by $( E, F, \pi, X )$,
where $( E, \pi, X )$ is a bundle, while $F$ is a topological space, such that \\

(3.2)~~~ $ F ~\mbox{and}~ E_x ~~\mbox{are homeomorphic, for}~ x \in X $ \\

and each $x \in X$ has an open neighbourhood $U \subseteq X$, such that \\

(3.3)~~~ $ U \times F ~~\mbox{and}~ \pi^{-1} ( U ) \subseteq E ~~\mbox{are homeomorphic} $ \\

{\bf Example 3.1.} \\

A trivial example, called the {\it trivial fiber bundle}, is given by \\

(3.4)~~~ $ E = X \times F,~~~ \pi = pr_X $

\hfill $\Box$ \\

In view of that, it follows that a fiber bundle is a bundle which
locally is a trivial fiber
bundle. \\

Here again, we can see that the concept of fiber bundle is very much
different from both that
of manifold, and of sheaf. \\
Regarding the difference with sheaf, we note that in the case of a
fiber bundle, the bundle or total space $E$ is homeomorphic not with
the base space $X$, but with the product $X \times
F$. \\

For both sheaves $( {\cal S}, \pi, X )$ and bundles $( E, \pi, X )$,
the respective spaces
${\cal S}$ and $E$ are larger than $X$. However, they are larger in very different ways. \\
Indeed, for a sheaf $( {\cal S}, \pi, X )$, even if ${\cal S}$ is
larger than $X$, it is
nevertheless locally homeomorphic with it. \\
On the other hand, for a fiber bundle  $( E, F, \pi, X )$, the space
$E$ is even locally much larger than $X$, since it is locally
homeomorphic not with $X$, but with $X \times F$, see
(3.3) above. \\

{\bf Example 3.2.} \\

A basic and useful nontrivial example is given by the {\it tangent
bundle} of a manifold. For an intuitively easy illustration, let us
consider it in the simple case of the $2$-manifold
$S^2 \subset \mathbb{R}^3$. \\
In this case, we associate to each $x \in S^2$ its tangent plane
$T_x$, and then take as a
fiber bundle \\

(3.5)~~~ $ E = \bigcup_{x \in S^2} \{ x \} \times T_x,~~~ F = \mathbb{R}^2,~~~ X = S^2 $ \\

and \\

(3.6)~~~ $ \pi : E \ni ( x, v ) \longmapsto x \in S^2 $ \\

The clearly \\

(3.7)~~~ $ E_x = \pi^{-1} ( x ) = \{ x \} \times T_x,~~~ x \in X = S^2 $ \\

hence (3.2) holds. \\

{\bf Example 3.3.} \\

As a further example, we show how the Moebius Band can be
represented as a fiber bundle. We
take \\

(3.8)~~~ $ X = S^1,~~~ F = [ 0, 1 ] $ \\

However, it is obvious that we cannot take $E = X \times F = S^1
\times [ 0, 1 ]$, since such
a choice would give a cylinder, and not the Moebius Band. \\
Therefore, we take $E$ in the following manner. We first undo the
circle $S^1$ and turn it into the interval $I = [ 0, 2 \pi ]$. Then
we consider the auxiliary topological space $D = I
\times F$ and define on it the following equivalence relation $\approx$ \\

(3.9)~~~ $ ( y, v ) \approx ( y\,', v\,' ) ~~\Longleftrightarrow~~
                  \left ( ~ \begin{array}{l}
                                 1)~~ y = y\,',~~~ v = v\,' \\ \\
                                 \mbox{or} \\ \\
                                 2)~~ y = 0,~~ y\,' = 2 \pi,~~ v + v\,' = 1 \\ \\
                                 \mbox{or} \\ \\
                                 3)~~ y = 2 \pi,~~ y\,' = 0,~~ v + v\,' = 1 \\ \\
                              \end{array} ~ \right ) $ \\

Finally, we take $E$ as the topological quotient space \\

(3.10)~~~ $ E = D / \approx $

\hfill $\Box$ \\

As seen in the case of the Moebius Band, the concept of fiber bundle
can easily deal with "curvature" or "non-flat" spaces, and does so
in a far more simple manner than the concept of
manifold. \\
However, both concepts have their comparative advantages. For
instance, manifolds can deal with spaces which are not representable
conveniently as fiber bundles. One such simple example is the
$2$-dimensional manifold given by the surface of a $3$-dimensional
ball with two holes
which penetrate it completely. \\ \\

{\bf 4. SINGULARITIES} \\

A major advantage of presheaves, when compared with manifolds and
bundles, is the surprisingly easy and general manner they can deal
with very large classes of {\it singularities}, and do
so without any analytic means. \\
In this regard, the essential role is played by the concept of {\it
flabby sheaf} in the original terminology of Jean Leray, or
equivalently, of {\it flabby presheaf}, in the present
terminology. \\

{\bf Definition 4.1.  FLABBY PRESHEAVES} \\

A presheaf $( S ( U ), \rho^U_V )$ on a topological space $X$ is
called {\it flabby}, if and
only if, for every open subsets $V \subseteq U \subseteq X$, we have \\

(4.1)~~~ $ \rho^U_V : S ( U ) \longrightarrow S ( V ) ~~\mbox{is surjective} $ \\

This rather unimpressive condition proves, among others, to go very
deep into the issue of
singularities, as illustrated in the following examples. \\

{\bf Example 4.1.} \\

Let $X = \mathbb{R}$ and $S ( U ) = {\cal C}^l ( X, \mathbb{R} )$,
for open $U \subseteq X$, where $0 \leq l \leq \infty$ is given.
Further, let $\rho^U_V$ be the usual restriction of ${\cal
C}^l$-smooth functions from the larger open subset $U$ to the
smaller one $V$. Then it is easy to see that, whenever we have open
subsets $V \subsetneqq U \subseteq X$, we shall
have \\

(4.2)~~~ $ \rho^U_V : S ( U ) \longrightarrow S ( V ) ~~\mbox{is not surjective} $ \\

Indeed, let us take for instance $V = ( 0, 1 )$ and $U = ( -1, 1 )$. Then for the function \\

(4.3)~~~ $ f ( x ) = 1 / x, ~~~ x \in V = ( 0, 1 ) $ \\

we obviously have $f \in {\cal C}^l ( V )$, for every $0 \leq l \leq
\infty$. However, there is no function $g \in {\cal C}^l ( U )$, for
no matter which $0 \leq l \leq \infty$, such that
$f = \rho^U_V ( g )$. \\

The reason for the fact that (4.1) does not hold for any $0 \leq l
\leq \infty$, is obviously
in the singularity of $f$ near to $0 \in X = \mathbb{R}$. \\

However, this has nothing to do with the unboundedness of $f$ near
to $0 \in X = \mathbb{R}$.
Indeed, the same situation occurs for the function \\

(4.4)~~~ $ h ( x ) = \sin ( 1 / x ), ~~~ x \in V = ( 0, 1 ) $ \\

Now contrary to presheaves which are not flabby, the flabby
presheaves can easily deal with a
very large class of singularities, as seen next. \\

{\bf Example 4.2.} \\

It has recently been shown that, given any $0 \leq l \leq \infty$,
then the smallest flabby presheaf which contains the non-flabby
presheaf $( {\cal C}^l ( U ), \rho^U_V )$, where $X =
\mathbb{R}^n$, is given by \\

(4.5)~~~ $ ( {\cal C}^l_{nd} ( U ), \rho^U_V ) $ \\

where \\

(4.6)~~~ $ {\cal C}^l_{nd} ( U ) = \left \{~ f : U \longrightarrow
\mathbb{R} ~~
                                                 \begin{array}{|l}
                        \exists~~ \Gamma \subset U,~~\mbox{closed, nowhere dense} ~: \\ \\
                        ~~~ f|_{U \setminus \Gamma} \in {\cal C}^l ( U \setminus \Gamma )
                                                  \end{array} ~ \right \} $ \\ \\

Obviously we have the strict inclusions \\

(4.7)~~~ $ {\cal C}^l ( U ) \subsetneqq {\cal C}^l_{nd} ( U ) $ \\

for $0 \leq l \leq \infty$ and open subsets $U \subseteq \mathbb{R}^n$. \\

What is important to note here is how much larger are the sets of
the functions in the right hand term of (4.7), than those in the
left hand term. In other words, how large is the class
of singular functions in ${\cal C}^l_{nd} ( U )$. \\
Indeed, there are two aspects involved here which make the sets of
functions ${\cal C}^l_{nd}
( U )$ considerably larger than the sets of functions ${\cal C}^l ( U )$. \\

First, a function $f \in {\cal C}^l_{nd} ( U )$, for which $f \in
{\cal C}^l ( U \setminus \Gamma )$, has its singularity set the
whole of $\Gamma$. And it is well known that the Lebesgue measure of
a closed nowhere dense subsets $\Gamma \subset U$ can be arbitrarily
near
to the Lebesgue measure of $U$. \\

Second, the condition $f \in {\cal C}^l ( U \setminus \Gamma )$ does
not impose any restrictions whatsoever on the function $f$ with
respect to its behaviour near to $\Gamma$.

\hfill $\Box$ \\

We turn now to integrable functions. \\

{\bf Example 4.3.} \\

Let $X = \mathbb{R}^n$ and $S ( U ) = {\cal L}^p ( U )$, for open $U
\subseteq X$, where $1 \leq p \leq \infty$ is given. Further, let
$\rho^U_V$ be the usual restriction of ${\cal L}^p$ functions from
the larger open subset $U$ to the smaller one $V$. Then the presheaf
$( S
( U ), \rho^U_V )$ is flabby. \\

Indeed, given $f \in {\cal L}^p ( V )$, we define $g : U
\longrightarrow \mathbb{R}$ simply by
$g|_V = f$ and $g = 0$ on $U \setminus V$. Then obviously $g \in {\cal L}^p ( U )$. \\

On the other hand, with local integrability, we have \\

{\bf Example 4.4.} \\

Let $X = \mathbb{R}^n$ and $S ( U ) = {\cal L}^p_{loc} ( U )$, for
open $U \subseteq X$, where $1 \leq p \leq \infty$ is given.
Further, let $\rho^U_V$ be the usual restriction of ${\cal
L}^p_{loc}$ functions from the larger open subset $U$ to the smaller
one $V$. Then the
presheaf $( S ( U ), \rho^U_V )$ is not flabby. \\

Indeed, with $X = \mathbb{R}$, let us take for instance $V = ( 0, 1
)$ and $U = ( -1, 1 )$.
Then for the function \\

(4.8)~~~ $ f ( x ) = 1 / x, ~~~ x \in V = ( 0, 1 ) $ \\

we obviously have $f \in {\cal L}^l_{loc} ( V ) \setminus {\cal L}^l
( V )$. Thus there cannot exist $g \notin {\cal L}^l_{loc} ( U )$,
such that $f = g|_V$, since $f$ is not integrable on any interval $(
0, a )$, with no matter how small $a > 0$.

\hfill $\Box$ \\

The above four examples are instructive with respect to a better
understanding of the essence of the flabbiness property. Indeed, in
Examples 4.1. and 4.2., we have the strict
inclusions \\

(4.9)~~~ $ {\cal C}^l ( U ) \subsetneqq {\cal C}^l_{nd} ( U ) $ \\

for $0 \leq l \leq \infty$ and open subsets $U \subseteq
\mathbb{R}^n$. And the smaller sets
of functions in (4.8) do not form flabby presheaves, while the larger sets of functions do. \\

On the other hand, we have the strict inclusions \\

(4.10)~~~ $ {\cal L}^p ( U ) \subsetneqq {\cal L}^p_{loc} ( U ) $ \\

for all $1 \leq p \leq \infty$ and open subsets $U \subseteq
\mathbb{R}^n$. And in this case, it is the smaller sets of functions
which form flabby presheaves, while the larger ones do
not. \\

Also, it should be noted that the functions in the presheaves $(
{\cal C}^l ( U ), \rho^U_V )$
in (4.9) which do not form flabby presheaves are defined by local conditions. \\

On the contrary, the functions in the presheaves $( {\cal L}^p ( U
), \rho^U_V )$ in
(4.10) which form flabby presheaves are defined by global conditions. \\

As had been noted elsewhere, the difference between presheaves with
functions defined by local, or on the contrary, global conditions,
can be significant. For instance, presheaves of
the first kind tend to be complete, while those of the second kind tend to fail to be so. \\

{\bf Open Problem.} \\

Similar to the extension from ${\cal C}^l$ to ${\cal C}^l_{nd}$ in
Examples 4.1., 4.2., let us
define for open $U \subseteq X$ and $1 \leq p \leq \infty$, the set of functions \\

(4.11)~~~ $ {\cal L}^p_{nd, loc} ( U ) = \left \{~ f : U
\longrightarrow \mathbb{R} ~~
                                                 \begin{array}{|l}
                        ~\exists~~ \Gamma \subset U,~~\mbox{closed, nowhere dense} ~: \\ \\
                        ~~~~ f|_{U \setminus \Gamma} \in {\cal L}^p_{loc} ( U \setminus \Gamma )
                                                  \end{array} ~ \right \} $ \\ \\

Then obviously \\

(4.12)~~~ $ {\cal L}^p ( U ) \subsetneqq {\cal L}^p_{loc} ( U ) \subsetneqq {\cal L}^p_{nd, loc} ( U ) $ \\

and the following two questions arise \\

1)~~ Is $ ( {\cal L}^p_{nd, loc} ( U ), \rho^U_V ) $ a flabby presheaf ? \\

2)~~ And in case it is, is it also the smallest flabby presheaf containing \\
     \hspace*{0.65cm} $ ( {\cal L}^p_{loc} ( U ), \rho^U_V ) $ ? \\

\newpage

{\bf APPENDIX : DIRECT LIMITS OF SETS \\
                \hspace*{2.9cm} AND OF MAPPINGS} \\

Let $( I, \leq )$ be a right directed preorder. A corresponding {\it right directed system of
sets} is given by \\

(A.1)~~~ $ ( E_\alpha, f_{\beta,\alpha} ),~~~ \alpha, \beta \in I,~~ \alpha \leq \beta $ \\

such that \\

(A.2)~~~ $ f_{\alpha,\alpha} = id_{E_\alpha},~~~ \alpha \in I $ \\

(A.3)~~~ $ f_{\gamma,\beta} \circ f_{\beta,\alpha} = f_{\gamma,\alpha},~~~
                          \alpha, \beta, \gamma \in I,~~ \alpha \leq \beta\leq \gamma $ \\

Assuming the sets $E_\alpha$ pair-wise disjoint, let \\

(A.4)~~~ $ G = \sum_{\alpha\in I} E_\alpha $ \\

and define on $G$ the equivalence relation $\approx$, by \\

(A.5)~~~ $ x \approx y ~~~\Longleftrightarrow~~~
                       \left (~ \begin{array}{l}
                                   \exists~~ \alpha, \beta, \gamma \in I ~:~ \\ \\
                                    ~~~~~~ *)~~ \alpha \leq \gamma,~~ \beta \leq \gamma \\ \\
                                    ~~~~ **)~~ x \in E_\alpha,~~ y \in E_\beta \\ \\
                                    ~~ ***)~~ f_{\gamma,\alpha} ( x ) = f_{\gamma,\beta} ( y )
                                 \end{array} ~ \right ) $ \\ \\

Then we define \\

(A.6)~~~ $ E = \underrightarrow{\lim}_{\,\alpha \in I\,} E_\alpha = G / \approx $ \\

Let \\

(A.7)~~~ $ f : G \longrightarrow E = G / \approx $ \\

be the canonical quotient mapping. For $\alpha \in I$, we define the {\it canonical direct
limit mapping} \\

(A.8)~~~ $ f_\alpha : E_\alpha \longrightarrow E = G / \approx $ \\

simply by the restriction \\

(A.9)~~~ $ f_\alpha = f|_{E_\alpha} $ \\

It follows easily that \\

(A.10)~~~ $ f_\beta \circ f_{\beta,\alpha} = f_\alpha,~~~
                    \alpha, \beta \in I,~~ \alpha \leq \beta $ \\

Let now be given any set $F$, and a family of mappings \\

(A.11)~~~ $ u_\alpha : E_\alpha \longrightarrow F $ \\

such that \\

(A.12)~~~ $ u_\beta \circ f_{\beta,\alpha} = u_\alpha,~~~
                               \alpha, \beta \in I,~~ \alpha \leq \beta $ \\

Then \\

(A.13)~~~ $ \begin{array}{l}
                  \exists~!~~ u : E \longrightarrow F ~: \\ \\
                  ~~~~~~ u_\alpha = u \circ f_\alpha,~~~ \alpha \in I
             \end{array} $ \\ \\

(A.14)~~~ $ u ~~\mbox{surjective} ~~~\Longleftrightarrow~~~
                       F = \bigcup_{\alpha \in I} u_\alpha ( E_\alpha ) $ \\ \\

(A.15)~~~ $ u ~~\mbox{injective} ~\Longleftrightarrow~
                    \left ( \begin{array}{l}
                                \forall~~ \alpha \in I,~~ x, y \in E_\alpha ~: \\ \\
                                ~~~~ u_\alpha ( x ) u_\alpha ( y ) ~\Longrightarrow~
                                  \left (~ \begin{array}{l}
                                              \exists~~ \beta \in I ~: \\ \\
                                              ~~~~ *)~~ \alpha \leq \beta \\ \\
                                              ~~~ **)~~ f_{\beta,\alpha} ( x ) =
                                                            f_{\beta,\alpha} ( y )
                                            \end{array} \right )
                              \end{array} \right ) $ \\ \\

{\bf Proof of (A.13)} \\

Let $v : G \longrightarrow F$ be such that $v_{E_\alpha} = u_\alpha$, with $\alpha \in I$.
Then (A.12) means that $v$ is compatible with the equivalence relation $\approx$ on $G$, hence
(A.13). \\


\begin{thebibliography} {99}

\bibitem{} Mallios A : Geometry of Vector Sheaves: An Axiomatic Approach to Differential
Geometry, Volumes 1,2. Springer-Verlag, New York, 1998

\end{thebibliography}
\end{document}